\documentclass[11pt,twoside]{article}
\usepackage{amssymb}
\usepackage{graphicx}

\usepackage{graphicx,type1cm,eso-pic,color}
\usepackage{amssymb}
\usepackage{amssymb,amsmath}
\usepackage{graphicx,epsfig}
\usepackage{enumerate}
\usepackage{color}
\usepackage {multicol}

\numberwithin{equation}{section}
\newtheorem{theorem}{Theorem}
\newtheorem{lemma}{Lemma}

\newtheorem{definition}{Definition}
\newtheorem{proposition}{Proposition}
\newtheorem{remark}{Remark}
\newtheorem{corollary}{Corollary}
\newtheorem{example}{Example}

\makeatletter
\AddToShipoutPicture{
	\setlength{\@tempdimb}{.5\paperwidth}
	\setlength{\@tempdimc}{.5\paperheight}
	\setlength{\unitlength}{1pt}
	\put(\strip@pt\@tempdimb,\strip@pt\@tempdimc)
}
\makeatother

\usepackage{scalerel,stackengine}
\stackMath
\newcommand\reallywidehat[1]{%
	\savestack{\tmpbox}{\stretchto{%
			\scaleto{%
				\scalerel*[\widthof{\ensuremath{#1}}]{\kern-.6pt\bigwedge\kern-.6pt}%
				{\rule[-\textheight/2]{1ex}{\textheight}}
			}{\textheight}%
		}{0.5ex}}%
	\stackon[1pt]{#1}{\tmpbox}%
}
\begin{document}
	\setcounter{page}{1}

	\thispagestyle{empty}
	\markboth{}{}

	\pagestyle{myheadings}
	\markboth{On weighted cumulative residual extropy and ...........}{ N. Gupta et al.}
	
	\date{}
	
	
	\noindent  
	
	\vspace{.1in}
	
	{\baselineskip 20truept
		
		\begin{center}
			{\Large {\bf On weighted cumulative residual extropy and weighted negative cumulative extropy }} \footnote{\noindent	{\bf *} E-mail: nitin.gupta@maths.iitkgp.ac.in\\
				{\bf ** } Corresponding author E-mail: skchaudhary1994@kgpian.iitkgp.ac.in\\
				{\bf *** } E-mail: pradeep.maths@kgpian.iitkgp.ac.in }\\
			
	\end{center}}

	\vspace{.1in}
	
	\begin{center}
		{\large {\bf Nitin Gupta*, Santosh Kumar Chaudhary**, Pradeep Kumar Sahu***}}\\
		{\large {\it Department of Mathematics, Indian Institute of Technology Kharagpur, West Bengal 721302, India }}
		\\
	\end{center}

	\vspace{.1in}
	\baselineskip 12truept

	\begin{abstract}
		In this paper, we define general weighted cumulative residual extropy (GWCRJ) and general weighted negative cumulative extropy (GWNCJ). We obtain its simple estimators for complete and right censored data. We obtain some results on GWCREJ and GWNCJ. We establish its connection to reliability theory and coherent systems. We also propose empirical estimators of weighted negative cumulative extropy (WNCJ).
		
		\vspace{.1in}
		
		
	\end{abstract}
	
	\section{Introduction}
	Entropy was introduced by Shannon \cite{shanon1948} as a measure of uncertainty involved in a random experiment. The entropy of a discrete random variable $X$ with probability mass function (pmf) $p_i,\text{for}\ \ i=1,2,\dots, n$ is defined as 
	\[ H(X)=-\sum_{i=0}^{n} p_i \log(p_i) .\]
	The entropy of $X$ for a non-negative continuous random variable $X$ with probability density function $f(.)$ is defined as 
	\[ H(X)=-\int_{0}^{\infty} f(x) \log(f(x))dx. \]
	Rao et al. \cite{raoetal2004}  introduced cumulative residual entropy (CRE) and studied its properties. CRE of a non-negative absolutely continuous random variable $X$ is defined as 
	\[ CRE(X)= -\int_{0}^{\infty} \bar{F}(x) \log(\bar{F}(x))dx.\]
	Extropy, as a dual measure of entropy, was introduced by Lad et al. \cite{lad2015}. The Extropy of a discrete random variable $X$ is defined as 
	\[ J(X)=-\sum_{i=0}^{n} (1-p_i) \log(1-p_i) .\]
	The Extropy of a non-negative continuous random variable $X$ is defined as 
	\[ J(X)=-\frac{1}{2} \int_{0}^{\infty} f^2(x)dx. \]
	Weighted extropy was introduced and studied independently by Bansal and Gupta  \cite{bansalgupta2020} and Balakrishnan et al. \cite{balaetal2020we}.
	Gupta and Chaudhary  \cite{guptaskc22} introduced and studied general weighted extropy. 
	Jahanshahi et al. \cite{janansahietal20} introduced cumulative residual extropy (CRJ) and studied its properties. CRJ of a non-negative absolutely continuous random variable $X$ is defined as 
	\[ CRJ(X)=-\frac{1}{2} \int_{0}^{\infty} \bar{F}^2(x)dx \]
	and cumulative past extropy (CPJ) for a non-negative absolutely continuous random variable is defined as 
	\[ CPJ(X)=-\frac{1}{2} \int_{0}^{\infty} F^2(x)dx. \]
	Kattumannil et al. \cite{kattumani22npecre} obtained a simple estimator of cumulative residual extropy for complete and right-censored data and studied their properties. In this paper, we also obtained a simple estimator of weighted cumulative residual extropy for complete and right censored data and studied their properties.
	Tahmasebi and Toomaj \cite{tahmasebitoomaj2020} introduced negative cumulative extropy (NCJ) and studied several properties of this concept. NCJ of a non-negative absolutely continuous random variable is defined as 
	\begin{equation}\label{ncextropydef}
		NCJ(X)=\bar \eta J(X)=\frac{1}{2} \int_{0}^{\infty} (1-F^2(x))dx. 
	\end{equation}
	For more properties of J(X ) see Lad et al. \cite{lad2015}. Recently, several researchers studied various forms of extropy and proposed some new measures. The Extropy of order statistics and record values is studied by Qiu \cite{qiu2017} and Jose and Sathar \cite{josesathar19} among others. Qiu and Jia (2018a) \cite{qiujia2018a} studied the residual extropy of order statistics. Qiu and Jia  \cite{qiujia2018b} developed the goodness of fit of uniform distribution using extropy. Qiu et al.  \cite{qiuetal2019} studied the properties of extropy for mixed systems. Jahanshahi et al. \cite{janansahietal20}   and Tahmasebi and Toomaj  \cite{tahmasebitoomaj2020} studied cumulative residual extropy and negative cumulative	extropy in detail. Recently, Balakrishnan et al.  \cite{balaetal2020we}, Bansal and Gupta  \cite{bansalgupta2020} and Sathar and Nair \cite{satharnair21a}, \cite{satharnair21b}, \cite{satharnair21c} studied different weighted versions of extropy. Kazemi et al. \cite{kazemietal.2022} introduced weighted cumulative past extropy and developed some characterization results. Several applications of extropy and its generalizations, such as in information theory, economics, communication theory, and physics, can be found in the literature. Here, we cite some references. Tahmasebi and Toomaj  \cite{tahmasebitoomaj2020} studied the stock market in OECD countries based on a generalization of extropy known as negative cumulative extropy. Balakrishnan et al. \cite{balakrishananetal2022} applied another version of extropy known as the Tsallis extropy to a pattern recognition problem. Kazemi et al.  \cite{kazemietal2021} explored an application of a generalization of extropy known as the fractional Deng extropy to a problem of classification. Tahmasebi et al. \cite{tahmasebietal2022} used some extropy measures for the problem of compressive sensing. The rest of this paper is organized as follows. Section \ref{s2gwcre} discusses general weighted cumulative residual extropy and its simple non-parametric estimator for a complete and right censored case. In Section \ref{s3gwce}, we define general weighted cumulative extropy and obtain a simple non-parametric estimator for a complete and right censored case. Section \ref{s4wnccpe} defines Weighted Negative cumulative past extropy and proposes some results. Section \ref{s5someinequality} is devoted to some inequalities which provide bounds for weighted negative cumulative past extropy. In Section \ref{s6ctoreth}, we established some connection of weighted negative cumulative extropy to reliability theory. Conditional weighted negative cumulative extropy is defined and studied in section \ref{s7cwpex}. Section \ref{s8empiricalofnwcex} provides two empirical estimators to weighted negative cumulative extropy. In section \ref{s9wncexofcs} we studied the weighted negative cumulative extropy of a coherent system. In the end, section \ref{s10conclusion} concludes this paper. 
	
	\section { General Weighted cumulative residual extropy} \label{s2gwcre}
	Let $X$ be a non-negative absolutely continuous random variable having density $f(\cdot)$, distribution function $F(\cdot)$ and survival function $\bar F(\cdot)=1-F(\cdot)$.
	\begin{definition}
		The general weighted cumulative residual extropy with weight $w(\cdot)\geq 0$ is given by 
		\[\eta^w J(X)=\frac{-1}{2}\int_{0}^{\infty}w(x) \bar F^2(x)dx\]
		If $w(x)=x^m$, then we represent $\eta^w J(X)$  as $\eta^m J(X)$.
	\end{definition}
	
	\begin{example}
		Let $X$ be an exponential random variable with a density function 
		$f_X(x)=\lambda e^{-\lambda x}, \ x> 0 , \ \lambda> 0$
		then for $m >-1,$ \[\eta^m J(X)=\frac{-1}{2}\int_{0}^{\infty} x^m \bar F^2(x)dx  =\frac{-1}{2}\int_{0}^{\infty} x^m e^{-2 \lambda x}dx =-\frac{\Gamma(m+1)}{2^{m+2} \lambda^{m+1}}.\]				
	\end{example}
	\subsection{Uncensored case}
	Let $X_1$ and $X_2$ be independent and identically distributed (i.i.d) random variables with common distribution function $F(\cdot)$. The survival function of $\min (X_1,X_2)$ is given by $\left( \bar F(x) \right)^2$. For the non-negative random variable $X$, we have
	\[\int_{0}^{\infty}x^mf(x)dx=m\int_{0}^{\infty}x^{m-1}\bar F(x)dx.\]
	Therefore we have
	\[\eta^m J(X)=\frac{-1}{2}\int_{0}^{\infty}x^m \bar F^2(x)dx=\frac{-1}{2(m+1)}E\left(\left(\min (X_1,X_2)\right)^{m+1}\right).\]
	Based on U-statistics an estimator of $\eta^m J(X)$ is given by
	\[T_{1,m}=\frac{-1}{n(n-1)(m+1)}\sum_{i=1}^{n-1}\sum_{j=i+1}^{n}\left(\min (X_i,X_j)\right)^{m+1}.\]
	Here $T_{1,m}$ is an unbiased estimator of $\eta^m J(X)$. Note that using order statistics we have
	\begin{equation}
		T_{1,m}=\frac{1}{n(n-1)(m+1)}\sum_{i=1}^{n-1}(i-n)X_{(i)}^{m+1},
	\end{equation}
	where $X_{(i)}$ be the $i$th order statistics from a random sample $X_1,\ldots,X_n$ from $F(\cdot)$. Now we study the asymptotic properties of $T_{1,m}$ using the limit theorems of U-statistics. $T_{1,m}$ is a consistent estimator of $\eta^m J(X)$, as $T_{1,m}$ is a U-statistic (Lehmann \cite{lehman1951}).
	\begin{theorem}\label{thmwcre}
		As $n\rightarrow \infty$, $\sqrt{n} (T_{1,m}-\eta^m J(X))$ is asymptotic normal with mean zero and variance $4\sigma_1^2$, where
		\begin{equation}\label{eqsigma1}
			\sigma_1^2=Var\left(X^{m+1}\bar F(X)+\int_{0}^{X}ydF(y)\right).
		\end{equation}
	\end{theorem}
	{\bf Proof.} Asymptotic normality follows from the central limit theorem of U-statistics (see Theorem 1, p. 76,  Lee \cite{lee2019}). Also, the asymptotic variance is $4\sigma_1^2$, where
	\[\sigma_1^2=Var\left(E\left(\left(\min(X_1,X_2)\right)^{m+1}|X_1\right)\right)\]
	(see Lee  \cite{lee2019}). Note that
	\begin{align*}
		E\left(\left(\min(X_1,X_2)\right)^{m+1}|X_1=x\right)&=E\left(x^{m+1}I(x<X_2)+X_2^{m+1}I(X_2\leq x)\right)\\
		&=x^{m+1}\bar F(x)+\int_{0}^{x}y^{m+1}dF(y).
	\end{align*}
	Hence the expression for variance is obtained as (\ref{eqsigma1}). Hence the result. \hfill $\blacksquare$
	
	\subsection{ Right censored case}\label{RCC1}
	Consider the case when the censoring times are independent of the lifetimes and the observations are randomly right-censored. Let K be the survival function of censoring random variable C. We want to find an estimator of the weighted cumulative residual extropy $\eta^m J(X)$ based on n independent and identical observations $\{(Y_i, \delta_i), 1 \leq i \leq n\},$ where $Y_i = \min(X_i, C_i)$ and $\delta_i = I(X_i \leq C_i),$ is the censoring indicator. A U-statistic defined for right censored data is given by 	\[T_{1c,m}=\frac{-1}{n(n-1)(m+1)}\sum_{i=1}^{n-1}\sum_{j<i;j=1}^{n} \frac{\left(\min (Y_i,Y_j)\right)^{m+1}\delta_i \delta_j}{\reallywidehat{K}(Y_i-) \reallywidehat{K}(Y_j-)},\]
	where $\reallywidehat{K}(.)$ is the Kaplan–Meier estimator of $K(.)$ ( see Datta et al. \cite{dataaetal2010}). The proof of the following theorem is similar to the proof of  Theorem 3 of Kattumannil et al. \cite{kattumannilsudheeshetal21}.
	\begin{theorem}
		As n $\rightarrow$ $\infty,\ \ T_{1c,m}$ converges in probability to $\eta^m J(X)$.
	\end{theorem}
	For deriving the asymptotic distribution of $T_{1c,m},$ let us define $N_i^c(t)=I(Y_i\leq t, \delta_i=0 )$ as the counting process corresponding to censoring random variable for the $i-$th subject and $R_i(u)=I(Y_i\geq u)$. Let $\lambda_c(t)$ be the hazard rate of the censoring variable $C$. The martingale associated with the counting process $N_i^c(t)$ is given by 
	\[M_i^c(t)=N_i^c(t)-\int_{0}^{t} R_i(u)\lambda_c(u)du.\]
	Let $H_c(x)=P(Y_1\leq x, \delta=1), y(t)=P(Y_1>t)$ and \[w(t)=\frac{1}{y(t)} \int \frac{h_1(x)}{\reallywidehat{K}(x-)}I(x>t)dH_c(x),\]
	where $h_1(y)=E(h(Y_1,Y_2)|Y_1=y).$ The proof of the following theorem follows from Datta et al. (2010) with choice of $h(Y_1, Y_2)={(\min(Y_1, Y_2))}^{m+1}.$
	
	\begin{theorem}
		Assume $E\left({(\min(Y_1,Y_2))}^{2(m+1)}\right) < \infty,\  \int \frac{h_1(x)}{\reallywidehat{K^2}(x-)} dH_c(x) <\infty$ and $\int_{0}^{\infty} w^2(t) \lambda_c(t)dt<\infty.$ As $n \rightarrow \infty,$ the distribution of $\sqrt{n} (T_{1c,m}-\eta^m J(X))$ is normal with mean zero and variance $4\sigma^2_{1c,m},$ where $\sigma^2_{1c,m}$ is given by \[\sigma^2_{1c,m}=Var\left(\frac{h_1(X)\delta_1}{\reallywidehat{K}(Y-)}+\int w(t)dM^c_1(t)\right)\]
	\end{theorem}

	\section{Weighted negative cumulative extropy}\label{s3gwce}
	
	In this section, we study the estimation of the general weighted negative cumulative extropy (GWNCJ). 
	\begin{definition}
		The general weighted negative cumulative extropy of a non-negative random variable $X$	is defined as 
		\[\bar \eta^wJ(X)=\frac{1}{2}\int_{0}^{\infty}w(x) \left(1- F^2(x)\right)dx\]
		If $w(x)=x^m$, then we represent $\bar \eta^wJ(X)$  as $\bar \eta^m J(X)$.
	\end{definition}

	\subsection{Uncensored case}
	Let $X_1$ and $X_2$ be independent and identically distributed random variables with common distribution function $F(\cdot)$. The distribution function of $\max (X_1,X_2)$ is given by 
	$\left(  F(x) \right)^2$. For the non-negative random variable $X$ we have
	\[\int_{0}^{\infty}x^mf(x)dx=m\int_{0}^{\infty}x^{m-1}\bar F(x)dx.\]
	Here
	\[\bar \eta^m J(X)=\frac{1}{2}\int_{0}^{\infty}x^m \left(1- F^2(x)\right)dx=\frac{1}{2(m+1)}E\left(\left(\max (X_1,X_2)\right)^{m+1}\right).\]
	Based on U-statistics an estimator of $\bar \eta^m J(X)$ is given by
	\[T_{2,m}=\frac{1}{n(n-1)(m+1)}\sum_{i=1}^{n-1}\sum_{j=i+1}^{n}\left(\max (X_i,X_j)\right)^{m+1}.\]
	Here $T_{2,m}$ is an unbiased and consistent estimator of $\bar \eta^m J(X)$. Note that using order statistics we have
	\begin{equation}
		T_{2,m}=\frac{1}{n(n-1)(m+1)}\sum_{i=1}^{n}(i-1) X_{(i)}^{m+1},
	\end{equation}
	Now we discuss the asymptotic distribution of  $T_{2,m}$ and the proof follows in a similar fashion as Theorem \ref{thmwcre}.
	\begin{theorem}\label{thmwce}
		As $n\rightarrow \infty$, $\sqrt{n} (T_{2,m}-\bar \eta^m J(X))$ is asymptotic normal with mean zero and variance $4\sigma_2^2$, where
		\begin{equation}\label{eqsigma2}
			\sigma_2^2=Var\left(X^{m+1} F(X)+\int_{X}^{\infty}y^{m+1}dF(y)\right).
		\end{equation}
	\end{theorem}
	
	\subsection{ Right censored case }\label{RCC2}
	Here, we find a simple estimator of $\bar \eta^m J(X)$ in the presence of right-censored observations. Using the same notations used in Section \ref{RCC1}, an estimator of  is given by
	
	\[T_{2c,m}=\frac{1}{n(n-1)(m+1)}\sum_{i=1}^{n-1}\sum_{j<i;j=1}^{n} \frac{\left(\max (Y_i,Y_j)\right)^{m+1}\delta_i \delta_j}{\reallywidehat{K}(Y_i-) \reallywidehat{K}(Y_j-)},\]
	where $\reallywidehat{K}(.)$ is the Kaplan–Meier estimator of $K(.).$
	\begin{theorem}
		As n $\rightarrow$ $\infty,\ \ T_{2c,m}$ converges in probability to $\bar \eta^m J(X).$
	\end{theorem}
	\begin{theorem}
		Let $h_1(y)=E\left((\max(Y_1,Y_2))^{m+1}| Y_1=y\right).$ Assume $E\left({(\max(Y_1,Y_2))}^{2(m+1)}\right) < \infty,\  \int \frac{h_1(x)}{\reallywidehat{K^2}(x-)} dH_c(x) <\infty$ and $\int_{0}^{\infty} w^2(t) \lambda_c(t)dt<\infty.$ As $n \rightarrow \infty,$ the distribution of $\sqrt{n} (T_{2c,m}-\bar{\eta}^m J(X))$ is normal with mean zero and variance $4\sigma^2_{2c,m},$  where $\sigma^2_{2c,m}$ is given by\[\sigma^2_{2c,m}=Var\left(\frac{h_1(X)\delta_1}{\reallywidehat{K}(Y-)}+\int w(t)dM^c_1(t)\right)\]
	\end{theorem}


	\section{Weighted Negative cumulative extropy} \label{s4wnccpe}
	Tahmasebi and Toomaj \cite{tahmasebitoomaj2020} defined negative cumulative extropy given in (\ref{ncextropydef}).
	In the last section, we defined general weighted negative cumulative extropy. Weighted Negative cumulative extropy (WNCJ)  is a particular case of general weighted Negative cumulative extropy when we choose $w(x)=x.$
	\begin{definition}
		Let $X$ be a non-negative absolutely continuous random variable with cdf $F$. We define the WNCJ of $X$ by
		\begin{equation}\label{wcpj}
			\bar\eta^1J(X)=\frac{1}{2}\int_{0}^{\infty}x \left(1-F^2(x)\right)dx.
		\end{equation}
	\end{definition}

	In this section, we study the properties of WNCPJ. Let us give a few examples before the main results.
	\begin{example}
		Let $X$ has $U[a,b]$ distribution. Then NCJ and WNCJ of the uniform distribution are
		
		\[\bar \eta J(X)=\frac{1}{2}\int_{0}^{\infty}\left(1-F^2(x)\right)dx=\frac{b-a}{6},\]
		and  \[\bar{\eta}^1J(X)=\frac{1}{2}\int_{0}^{\infty}x\left(1-F^2(x)\right)dx=\frac{b-a}{24}(5a+3b),\]	
		respectively.
	\end{example}
	
	\begin{example}
		Let $X$ have power distribution with pdf $f(x)=\lambda x^{\lambda-1},  x\in(0,1), \lambda > 1$.	
		Then $NCPJ$ and $WNCPJ$ of the distribution are 
		\begin{equation*}
			\bar \eta J(X)=\frac{1}{2(2\lambda+1)},\  \mbox{and}, \ \bar{\eta}^1J(X)=\frac{\lambda}{4(\lambda+1)}.
		\end{equation*}	
	\end{example}
	
	\begin{example}
		Let  $X$ have exponential distribution with mean $\frac{1}{\lambda}$ has the pdf $f(x)=\lambda e^{-\lambda x},  x>0, \lambda > 0$.
		Then $NCPJ$ and $WNCPJ$ of the distribution are
		\begin{equation}
			\bar \eta J(X)=\infty,\  \mbox{and}, \ \bar{\eta}^1J(X)=\frac{5}{4\lambda^{2}}.
		\end{equation}
	\end{example}
	Now we see the effect of linear transformation on WNCJ in the following proposition
	\begin{proposition}
		Let $X$ be a non-negative random variable. If $Y=aX+b,\ a>0,\ b\geq 0,$ then \[\bar \eta^1J(Y)=a^2 \bar\eta^1J(X)+ab \bar \eta J(X)\]
	\end{proposition}
	{\bf Proof: } The proof holds using (\ref{wcpj})
	and noting that $F_Y(y)=F_X\left(\frac{y-b}{a}\right),\ y>b$.\hfill $\blacksquare$
	
	Here we provide a lower bound for WNCJ in terms of extropy.
	\begin{theorem}\label{relation}
		Let $X$ be a random variable with pdf $f(\cdot)$ and extropy $J(X)$, then
		\begin{equation}
			\bar \eta^1  J(X)\geq C^* \exp\{2J(X)\},
		\end{equation} 
		where $C^*=\frac{1}{2}\exp \{E\left(\log\left(X\left(1-F^2(X)\right)\right)\right)\}$ and $\exp(x)=e^x.$
	\end{theorem}
	{\bf Proof: } Using the log-sum inequality, we have
	\begin{align*}
		\int_{0}^{\infty}f(x)\log\left(\frac{f(x)}{x(1-F^2(x))}\right) dx\geq -\log\left(\int_{0}^{\infty}x(1-F^2(x))dx\right).
	\end{align*}
	Then it follows that
	\begin{align*}
		\int_{0}^{\infty}f(x)\log f(x)dx-\int_{0}^{\infty}f(x)\log\left(x(1-F^2(x))\right) dx \geq -\log\left(\int_{0}^{\infty}x(1-F^2(x))dx\right).
	\end{align*}
	Note that $\log(f)< f$, hence
	\begin{align}\label{2.3aa}
		-\int_{0}^{\infty} f^2(x)dx+\int_{0}^{\infty}f(x)\log\left(x(1-F^2(x))\right) dx=&2J(X)+E\left(\log\left(X(1-F^2(X))\right)\right)\nonumber\\
		&\leq \log\left(2 \bar \eta^1J(X)\right).
	\end{align}
	Exponentiating both sides of (\ref{2.3aa}), we have
	\[\bar \eta^1J(X)\geq \frac{1}{2}exp\{2J(X)+E\left(\log\left(X(1-F^2(X))\right)\right)\}\]
	Hence the result.\hfill $\blacksquare$

	\begin{theorem}
		$\bar{\eta}^1J(X)=0$ if and only if $X$ is degenerate.
	\end{theorem}
	{\bf Proof} Suppose $X$ be degenerate at point $c$, then by using the definition of degenerate function and 
	$\bar{\eta}^1J(X)$, we have $\bar{\eta}^1J(X)=0$. Now consider $\bar{\eta}^1J(X)=0$, i.e., \[\int_{0}^{\infty}x(1-F^2(x))dx=0\].
	Noting that the integrand in the above integral is non-negative, we have $1-F^2(x)=0$, that is, $F(x)=1$ for almost all $x\in S$, where $S$ denote the support of 
	random variable $X$. Hence $X$ is degenerate.

	\section{Some inequalities}\label{s5someinequality}

	This section deals with obtaining the lower and upper bounds for WNCJ.
	\begin{remark}
		Consider $X$ to be a non-negative random variable. then
		\begin{equation}
			\bar \eta^1J(X)\geq \frac{E(X^2)}{4}
		\end{equation}
	\end{remark}
	{\bf Proof} Using inequality $F^2(x)\le F(x)$ and $E(X^2)=2\int_{0}^{\infty}x\bar F(x)dx$, the result follows.

	\begin{proposition}
		Consider a non-negative continuous random variable $X$ having cdf $F_X(\cdot)$ and support $[a,\infty), a>0$. Then
		\begin{equation}
			\bar \eta^1J(X)\geq a\bar \eta J(X).
		\end{equation}
	\end{proposition}
	{\bf Proof} Note that
	\begin{align*}
		\int_{a}^{\infty}x\left(1-F^2(x)\right)dx &\geq a\int_{a}^{\infty}\left(1-F^2(x)\right)dx\\
		\implies \ \bar \eta^1J(X) &\geq a\bar \eta J(X).
	\end{align*}   \hfill $\blacksquare$

	Consider two random variables $X$ and $Y$ having cdfs $F$ and $G$, respectively. Then $X \leq_{st}Y$ whenever $F(x)\geq G(x),\ \forall \ x\in \mathbb{R}$; where the notation $X \leq_{st}Y$ means that $X$ is less than or equal to $Y$ in usual stochastic order. One may refer to Shaked and Shanthikumar \cite{shakedshanti2007} for detail of stochastic ordering. In the following proposition, we show the ordering of WNCJ is implied by the usual stochastic order.
	\begin{proposition}
		Let $X_1$ and $X_2$ be non-negative continuous random variables. If $X_1\leq_{st}X_2$, then $ \bar \eta^1J(X_1)\leq  \bar \eta^1J(X_2)$.
	\end{proposition}
	{\bf Proof} Using $X_1\leq_{st}X_2$ and (\ref{wcpj}), the result follows.\hfill $\blacksquare$
	
	In the following, we obtain the WNCJ of the $n$-th order statistic. The WNCJ of the $n$-th order statistic is 
	\begin{equation}\label{largest}
		\bar \eta^1J(X_{n:n})=\frac{1}{2}\int_{0}^{\infty}x\left(1-F^{2}_{X_{n:n}}(x)\right)dx,
	\end{equation}
	where $F^{2}_{X_{n:n}}(x)=F^{2n}_{X}(x)$. Using transformation $u=F(x)$ in $(\ref{largest})$,
	\begin{equation}
		\bar \eta^1J(X_{n:n})= \frac{1}{2}\int_{0}^{1}\frac{\left(1-u^{2n}\right)F^{-1}(u)}{f(F^{-1}(u))} du,
	\end{equation}
	where $F^{-1}(x)$ is the inverse function of $F(x)$.
	\begin{example}
		Let $X$ have the uniform distribution on (0,1) with pdf $f(x)=1, \ x\in (0,1)$. Then $F^{-1}(u)=u,\ u\in (0,1)$ and $f(F^{-1}(u))=1, \ u\in (0,1)$, hence 
		$\bar \eta^1J(X_{n:n})=\dfrac{n}{4(n+1)}$.
	\end{example}
	\begin{example}
		let $X$ follow power-law distribution with pdf $f(x)=\lambda x^{\lambda-1}, \lambda > 1, x\in(0,1)$.Then $F^{-1}(u)=u^{\frac{1}{\lambda}},\ u\in (0,1)$ and $f(F^{-1}(u))= \lambda u^{\frac{\lambda-1}{\lambda}}, \ u\in (0,1)$, hence 
		$\bar \eta^1J(X_{n:n})=\dfrac{-n\lambda}{4(n\lambda +1)}$.
	\end{example}	
	\begin{remark}
		Consider $\Lambda =\bar \eta^1J(X_{n:n})-\bar \eta^1J(X)$. Since $F^{2n}(x)\leq F^{2}(x)$, hence $\Lambda \geq 0$.
	\end{remark}
	\section{Connection to reliability theory}\label{s6ctoreth}

	Consider a non-negative continuous random variable $X$ with cdf $F$, such that $E(X)<\infty$. The mean inactivity time (MIT) of $X$ is defined as
	\begin{equation}
		MIT(t)=\int_{0}^{t}\dfrac{F(x)}{F(t)}dx, \ t\geq 0.
	\end{equation}

	The MIT function finds many applications in reliability, forensic science, and so on. In the following theorem, we show the relationship between WNCJ and the second moment of inactivity time (SMIT) function. For detail about SMIT one may refer Kundu and Nanda \cite{kundunanda2010}.
	\begin{definition}
		Let $X$ be a non-negative continuous random variable. Then SMIT is 
		\begin{equation}\label{smit}
			SMIT(t)=E\left((t-X)^2|X\leq t\right)=2tMIT(t)-\int_{0}^{t}2x\dfrac{F(x)}{F(t)}dx, \ t\geq 0.
		\end{equation}
	\end{definition}

	\begin{theorem}
		Let $X$ be a non-negative random variable with finite WNCJ. Then we have
		\begin{equation*}
			\bar \eta^1J(X)= \frac{1}{2}\left[E(X^2)+E(XF(X) MIT(X))-\frac{1}{2}E(F(X)SMIT(X))\right].
		\end{equation*}
	\end{theorem}
	{\bf Proof}
	Consider
	\begin{align}\label{smit2}
		\bar \eta^1J(X)&=\frac{1}{2}\int_{0}^{\infty}x(1-F^2(x))dx\nonumber\\
		&=\frac{1}{2}\int_{0}^{\infty}x\bar F(x)(1+F(x))dx\nonumber\\
		&=\frac{1}{2}\left[\int_{0}^{\infty}x\bar F(x)dx+\int_{0}^{\infty}xF(x)\bar F(x)dx\right]\nonumber\\
		&=\frac{1}{2}\left[E(X^2)+\int_{0}^{\infty}xF(x)\left(\int_{x}^{\infty}f(t)dt\right)dx\right]\nonumber\\
		&=\frac{1}{2}\left[E(X^2)+\int_{0}^{\infty}f(t)\left(\int_{0}^{t}xF(x)dx\right)dt\right]
	\end{align}
	Now using (\ref{smit}) and (\ref{smit2}), we have
	\begin{align*}
		\bar \eta^1J(X)&=\frac{1}{2}\left[E(X^2)+\int_{0}^{\infty}\frac{f(t)F(t)}{2}\left\{2t MIT(t)-SMIT(t)\right\}dt\right]\\
		&=\frac{1}{2}\left[E(X^2)+E(XF(X) MIT(X))-\frac{1}{2}E(F(X)SMIT(X))\right]
	\end{align*}
	\hfill $\blacksquare$

	A bound for $\bar \eta^1J(X)$ can be provided in terms of the hazard rate function.
	\begin{proposition}
		Let $X$ be a non-negative continuous random variable with finite hazard rate function $h(\cdot)$ and $\bar \eta^1J(X)$. Then,
		\begin{equation}
			\bar \eta^1J(X)\geq E(S(X)),
		\end{equation}
		where $S(t)=\int_{0}^{t}x\left( \int_{0}^{x}h(v)dv\right)dx$.
	\end{proposition}
	{\bf Proof}
	From (\ref{smit2}) we have
	\begin{align*}
		\bar \eta^1J(X)&=\frac{1}{2}\left[E(X^2)+\int_{0}^{\infty}f(t)\left(\int_{0}^{t}xF(x)dx\right)dt\right]\\
		&\le \frac{1}{2}\left[E(X^2)-\int_{0}^{\infty}f(t)\left(\int_{0}^{t}x \log \bar F(x)dx\right)dt\right]\\
		&= \frac{1}{2}\left[E(X^2)-\int_{0}^{\infty}f(t)\left(\int_{0}^{t}x\left( \int_{0}^{x}h(v)dv\right)dx\right)dt\right]\\
		&=\frac{1}{2}\left[E(X^2)-E(S(X))\right]
	\end{align*}
	where  $S(t)=\int_{0}^{t}x\left( \int_{0}^{x}h(v)dv\right)dx$. Hence the result.\hfill $\blacksquare$

	\section{Conditional weighted negative cumulative extropy}
	\label{s7cwpex}
	
	Now we consider the conditional Negative weighted cumulative extropy (CWNCJ). Consider a random variable $Z$ on probability space $(\Omega, \mathbb{A}, P)$ such that $E|Z|<\infty$.
	The conditional expectation of $Z$ given sub $\sigma$-field $\mathbb{G}$, where $\mathbb{G}\subseteq \mathbb{A}$, is denoted by $E(Z|\mathbb{G})$. For the random variable $I_{(Z\leq z)}$, we denote $E(I_{(Z\leq z)}|\mathbb{G})$ by $F_Z(z|\mathbb{G})$.

	\begin{definition}
		For a non-negative random variable $X$, given $\sigma$-field $\mathbb{G}$, the CWNCJ $\bar \eta^1J(X|\mathbb{G})$ is defined as 
		\begin{align}
			\bar \eta^1J(X|\mathbb{G})=\frac{1}{2}\int_{0}^{\infty}x \left(1-F^2_X(x|\mathbb{G})\right)dx.
		\end{align}
	\end{definition}
	
	Now we assume that the random variables are continuous and non-negative.
	
	\begin{lemma}
		If $\mathbb{G}$ is a trivial $\sigma$-field, then $\bar \eta^1J(X|\mathbb{G})=\bar \eta^1J(X)$.
	\end{lemma}
	{\bf Proof} Since here $F_X(x|\mathbb{G})=F_X(x)$, then the proof follows.\hfill $\blacksquare$

	\begin{proposition}\label{propG}
		If $X\in L^p$ for some $p>2$, then $E[\bar \eta^1J(X|\mathbb{G})|\mathbb{G}^*]\leq \bar \eta^1J(X|\mathbb{G}^*)$, provided that $\mathbb{G}^*\subseteq \mathbb{G}$.
	\end{proposition}
	{\bf Proof}
	Consider
	\begin{align*}
		E[\bar \eta^1J(X|\mathbb{G})|\mathbb{G}^*]&=\frac{1}{2}\int_{0}^{\infty}x\left[1-E\left(\left[P(X\leq x|\mathbb{G})\right]^2|\mathbb{G}^*\right)\right]dx\\
		&\leq \frac{1}{2}\int_{0}^{\infty}x\left[1-\left[E\left(P(X\leq x|\mathbb{G})|\mathbb{G}^*\right)\right]^2\right]dx\\
		&= \frac{1}{2}\int_{0}^{\infty}x\left[1-\left[E\left(E(I_{(X\leq x)}|\mathbb{G})|\mathbb{G}^*\right)\right]^2\right]dx\\
		&= \frac{1}{2}\int_{0}^{\infty}x\left[1-\left[E\left(I_{(X\leq x)}|\mathbb{G}^*\right)\right]^2\right]dx\\
		&=\frac{1}{2}\int_{0}^{\infty}x\left[1-F^2_X(x|\mathbb{G}^*)\right]dx\\
		&=\bar \eta^1J(X|\mathbb{G}^*),
	\end{align*}
	where the second step follows using Jensen's inequality for convex function $\phi(x)=x^2$. Hence the result.\hfill $\blacksquare$
	
	In the following theorem, we investigate the relationship between conditional extropy and $\bar \eta^1J(X|\mathbb{G})$.
	\begin{theorem}
		Let $\bar \eta^1J(X|\mathbb{G})$ is conditional negative cumulative extropy. Then we have
		\begin{equation}
			\bar \eta^1J(X|\mathbb{G})\geq B^* exp\{2J(X|\mathbb{G})\},
		\end{equation} 
		where $B^*=\frac{1}{2}exp\{E\left(\log \left(X \left(1-F^2(X)\right)\right)|\mathbb{G}\right)\}$
	\end{theorem}
	{\bf Proof}
	The proof is on the similar lines as of Theorem \ref{relation}, hence omitted.\hfill $\blacksquare$

	\begin{theorem}
		For a random variable $X$ and $\sigma$-field $\mathbb{G}$, we have
		\begin{align}\label{thmprop}
			E\left(\bar \eta^1J(X|\mathbb{G})\right)\leq \bar \eta^1J(X),
		\end{align}
		and the equality holds if and only if $X$ is independent of $\mathbb{G}$.
	\end{theorem}
	{\bf Proof}
	If in Proposition \ref{propG}, $\mathbb{G}^*$ is trivial $\sigma$-field, then (\ref{thmprop}) can be easily obtained. Now assume that $X$ is independent of 
	$\mathbb{G}$, then 
	\begin{align}\label{eqthmprop1}F_X(x|\mathbb{G})&=F_X(x)\nonumber\\
		\implies\ \bar \eta^1J(X|\mathbb{G})&=\bar \eta^1J(X).
	\end{align}
	On taking expectation to both sides of  (\ref{eqthmprop1}), we get equality in  (\ref{thmprop}). Conversely, assume that equality in  (\ref{thmprop}) holds. It is sufficient to show that 
	$F_X(x|\mathbb{G})=F_X(x)$, to prove independence between $X$ and $\sigma$-field $\mathbb{G}$. Take $U=F_X(x|\mathbb{G})$, and since the function $\phi(u)=u^2$
	is convex hence $E(U^2)\geq E^2(U)=F_X^2(x)$, and  also due to equality in (\ref{thmprop}), we have
	\begin{align*}
		\int_{0}^{\infty}x\left(1-E(U^2)\right)dx=\int_{0}^{\infty}x\left(1-F^2_X(x)\right)dx=\int_{0}^{\infty}x\left(1-E^2(U)\right)dx.
	\end{align*}
	Hence $E(U^2)=E^2(U)$. Now using Corollary 8.1 of Hashempour et al. \cite{hashempouretal22}, we have $F_X(x|\mathbb{G})=F_X(x)$. Hence the proof.\hfill $\blacksquare$

	For the Markov property for non-negative random variables $X,\ Y$ and $Z$, we have the following proposition.
	\begin{proposition}
		Let $X\rightarrow Y \rightarrow Z$ is a Markov chain, then
		\begin{equation}\label{po1}
			\bar \eta^1J(Z|X,Y)=\bar \eta^1J(Z|Y)
		\end{equation}
		and
		\begin{equation}\label{po2}
			E\left(\bar \eta^1J(Z|Y)\right)\leq E\left(\bar \eta^1J(Z|X)\right).
		\end{equation}
	\end{proposition}
	{\bf Proof}
	By the definition of $\bar \eta^1J(Z|X,Y)$ and using the Markov property, (\ref{po1}) holds. \\
	Now letting $\mathbb{G}^*=\sigma(X),\ \mathbb{G}=\sigma(X,Y)$ and $X=Z$ in Proposition \ref{propG}, we have
	\begin{equation}\label{eq2}
		\bar \eta^1J(Z|X)\geq E\left(\bar \eta^1J(Z|X,Y)|X\right)
	\end{equation}
	Taking expecation on both sides of (\ref{eq2}), we have
	\begin{align*}
		E\left(\bar \eta^1J(Z|X)\right)&\geq E\left(E\left(\bar \eta^1J(Z|X,Y)|X\right)\right)\\
		&=E\left(\bar \eta^1J(Z|X,Y)\right)\\
		&=E\left(\bar \eta^1J(Z|Y)\right),
	\end{align*}
	where the last equality holds using (\ref{po1}). Hence the result (\ref{po2}) holds.\hfill $\blacksquare$

	
	\section{Empirical of WNCEX}\label{s8empiricalofnwcex}
	Here, we will find an estimator of weighted negative cumulative extropy by means of empirical weighted negative cumulative extropy. Let $X_1, X_2,\dots, X_n$ be a random sample from an absolutely continuous cumulative distribution function $F(x).$ If $X_{1:n}, X_{2:n}, \dots, X_{n:n}$ display the order statistics of random samples $X_1, X_2,\dots, X_n$, then the empirical measure of $F(x)$ for $i=1,2,\dots, n$ is given by
	$$	\reallywidehat{F}_n(X)=
	\begin{cases}
		0, \ x<X_{1:n}\\
		\frac{1}{n}, \ X_{i:n}\leq x < X_{i+1:n} \  \text{for} \ \  i=1,2,\dots, n-1\\
		1,\ x\geq X_{n:n} \ 
	\end{cases}
	$$
	The empirical measure of  $\bar \eta^1J(X)=\frac{1}{2} \int_{0}^{\infty} x(1-F^2(x))dx$  is obtained as
	\begin{align}
		\bar \eta_1^1J(\reallywidehat{F}_n) &=\frac{1}{2} \int_{0}^{\infty} x(1-\reallywidehat{F}_n^2(x))dx \nonumber\\
		&=\frac{1}{2}\sum_{i=1}^{n-1}\int_{X_{i:n}}^{X_{i+1:n}} \left[1-\left(\frac{i}{n}\right)^2\right] xdx \nonumber\\
		&=\frac{1}{4}\sum_{i=1}^{n-1} (X^2_{i+1:n}-X^2_{i:n}) \left[1-\left(\frac{i}{n}\right)^2\right]
	\end{align}
	Moreover, from equation (23) of Tahmasebi and Toomaj \cite{tahmasebitoomaj2020}, the empirical measure of negative cumulative extropy $\bar \eta J(X)= \frac{1}{2} \int_{0}^{\infty} (1-F^2(x))dx$ is
	\begin{align}
		\bar \eta J(\hat{F}_n)= \frac{1}{2}\sum_{i=1}^{n-1} (X_{i+1:n}-X_{i:n}) \left[1-\left(\frac{i}{n}\right)^2\right]
	\end{align}
	
	and from equation (16) of Jahanshahi et al. \cite{janansahietal20}, the empirical measure of negative cumulative residual extropy (NCRE) \[NCRE(X)= \frac{1}{2} \int_{0}^{\infty} \bar{F}^2(x)dx\] is 
	\begin{equation}
		NCRE(\hat{F}_n)= \frac{1}{2}\sum_{i=1}^{n-1} (X_{i+1:n}-X_{i:n}) \left[1-\frac{i}{n}\right]^2 
	\end{equation}
	We can write $\bar \eta^1J(X)$ as 
	\begin{align}
		\bar \eta^1J(X)&=\frac{1}{2} \int_{0}^{\infty} x(1-F^2(x))dx \nonumber \\
		&= \frac{1}{4} \int_{0}^{1} (1-u^2) \left[\frac{d(F^{-1}(u))^2}{du}\right] du
	\end{align}  
	The following theorem states that empirical weighted negative cumulative extropy $\bar \eta_1^1J(\reallywidehat{F}_n)$ converges almost surely to weighted negative cumulative extropy. Since almost sure convergence is a stronger condition than convergence in probability and convergence in distribution. Therefore $\bar \eta_1^1J(\reallywidehat{F}_n)$ is also a consistent estimator. Proof of the next statement follows from the proof of Theorem 4.1 of Tahmasebi and Toomaj  \cite{tahmasebitoomaj2020}.
	\begin{theorem}
		$\bar \eta_1^1J(\reallywidehat{F}_n)$ converges almost surely to $\bar \eta^1J(X)$.
	\end{theorem}
	
	One other way to find an estimator of $\bar \eta^1J(X)$ is to use the method suggested by Vasicek \cite{vasicek}. Following the idea of Vasicek  \cite{vasicek}, an estimator of $\bar \eta^1J(X)$ will be calculated by replacing the distribution function $F$ with an empirical distribution function $\reallywidehat{F}_N$ and using the difference operator in place of a differential operator. The derivative of $F^{-1} (u)$ with respect to $u$, that is, $\frac{dF^{-1} (u)}{du}$ will be estimated as 
	\begin{align*}
		\frac{X_{i+m:n}-X_{i-m:n}}{\reallywidehat{F}_n(X_{i+m:n})-\reallywidehat{F}_N(X_{i-m:n}}=\frac{X_{i+m:n}-X_{i-m:n}}{\frac{i+m}{n}-\frac{i-m}{n}}=\frac{X_{i+m:n}-X_{i-m:n}}{2m/n}.  
	\end{align*}
	Here window size $m$  is a positive integer less than $\frac{n}{2}$  and $X_{r:n}$ denotes $r$th order statistics from sample $X_1, X_2, \dots X_n.$  If $i+m>n$ then we consider $X_{i+m:n} = X_{n:n}$ and if $i+m<1$ then we consider $X_{i-m:n} = X_{1:n}.$ 
	
	Empirical estimator of $\bar \eta^1J(X)$ is obtained as 
	\begin{align}
		\bar \eta_2^1J(\reallywidehat{F}_n)= \frac{1}{4n}\sum_{i=1}^{n} \frac{(X^2_{i+m:n}-X^2_{i-m:n})}{2m/n}  \left[1-\left(\frac{i}{n+1}\right)^2\right]
	\end{align}
	The following theorem says empirical estimator $\bar \eta_2^1J (\reallywidehat{F}_n)$ is a consistent estimator of weighted negative cumulative extropy. Following lines of proof of Theorem 1 of Vasicek \cite{vasicek}, we obtain proof of the next theorem.
	\begin{theorem}
		$\bar \eta_2^1J(\reallywidehat{F}_n)$ converges in probability to $\bar \eta^1J(X)$.
	\end{theorem}


	\section{WNCEX of coherent system}\label{s9wncexofcs}
	A system with a monotone structure function and without any irrelevant components is called a coherent system ( for more on a coherent system, see Navarro \cite{navarro2022}, Barlo and Proschan \cite{barloproschan1975} ). Navarro et al. \cite{navaroetal2013} proved  \[F_T(x)=q(F_X(x))\] for a coherent system with identically distributed (i.d.) components where $F_T(x)$ is cdf of a coherent system, $F_X$ is the common cdf of the components and $q$ is distortion function. Function $q$ is an increasing continuous function in $[0,1]$ such that $q(0)=0$ and $q(1)=1$ which only depends on the structure of the system and the copula of the random vector $(X_1, X_2, \dots, X_n).$ From the weighted negative cumulative extropy defined in (\ref{wcpj}) and the probability integral transformation $V=F(T)$, we have
	\begin{align} \label{wcpecs}
		\bar \eta^1J(T)&=\frac{1}{2} \int_{0}^{\infty} x \left(1-F^2_T(x) \right)dx =  \int_{0}^{\infty} x \phi(F_T(x)) dx \nonumber\\
		&=\int_{0}^{\infty} x \phi(q(F(x))) dx= \int_{0}^{1} \frac{F^{-1}(v)\phi (q(v))}{f(F^{-1}(v))}dv,
	\end{align}
	where $\phi(u)=\frac{1-u^2}{2}, \ \text{for}\ \ 0<u<1.$ Also, we have
	\begin{align} \label{cpecs1}
		\bar \eta^1J(X)&=\frac{1}{2} \int_{0}^{\infty} x \left(1-F^2_X(x) \right)dx = \int_{0}^{1} \frac{F^{-1}(v)\phi (v)}{f(F^{-1}(v))}dv,
	\end{align}
	\begin{example}
		For a parallel system  with iid component of U(0,1) and  lifetime $T=\max(X_1, X_2, X_3,..., X_n)$,  we have $q(v)=v^n$, $\bar \eta^1J(T)=\frac{n}{4(n+1)}$ and $\bar \eta^1J(X)=\frac{1}{4}.$ Thus we note that $\bar \eta^1J(T) \leq \bar \eta^1J(X)$.
	\end{example}
	
	\begin{example}
		For a coherent system with identically distributed components of exponential distribution with cdf $F(x)=1-e^{-\frac{x}{\mu}}$ for $\mu > 0$ and $x>0$ and lifetime $T=\max\{\min\{X_1, X_2\}, \min\{X_3, X_4\}\},\ \ f(F^{-1}(v))=\frac{1-v}{\mu}.$ The maximal signature of the system is (0,4,-4,1) and so, from (\ref{wcpecs}) we obtain
		\begin{align}
			\bar \eta^1J(X)&= \int_{0}^{1} \frac{F^{-1}(v)\phi (q(v))}{f(F^{-1}(v))}dv \nonumber\\
			&= -\frac{\mu^2}{2} \int_{0}^{1} \frac{ln(1-v)\left[1-(4v^2-4v^3+v^4)^2\right]}{1-v} \nonumber \\
			&=0.3602 \mu^2.\nonumber
		\end{align}
	\end{example}
	
	\begin{proposition}
		Let T be the lifetime of the coherent system with i.d. components and distortion function $q$. If $\phi(q(u))\geq (\leq) \phi(u)$ for all $0<u<1,$ it holds that \[ \bar \eta^1J(T) \geq (\leq) \bar \eta^1J(X). \]
	\end{proposition}
	
	\begin{theorem}
		Let $T$ denote the lifetime of a coherent system with i.d. components with a common pdf $f$ and distortion function $q$. If $S$ is the support of $f$, $m=\inf_{x\in S} \frac{f(x)}{x} > 0$ and $M=\sup_{x\in S}\frac{f(x)}{x}<\infty,$ then 
		
		\[\frac{I_q}{M} \leq \bar \eta^1J(T) \leq \frac{I_q}{m}  \]
		where $I_q=\int_{0}^{1} \phi(q(u))du$ and $\phi(u)=\frac{1-u^2}{2} for 0<u<1.$
	\end{theorem}
	
	\begin{proposition}
		Assume that $T$ is the lifetime of a coherent system with i.d. components
		and with distortion function $q$. Then
		\[B_1 \bar\eta^1J(X_1) \leq \bar \eta^1J(T) \leq B_2 \bar \eta^1J(X_1)\]
		where $B_1=\inf_{u\in (0,1)} \left(\frac{\phi(q(u))}{\phi(u)}\right)$ and $B_2=\sup_{u\in (0,1)} \left(\frac{\phi(q(u))}{\phi(u)}\right).$		
	\end{proposition}
	
	\begin{proposition}
		Let $T_1$ and $T_2$ be the lifetimes of two coherent systems with i.d. components and with distortion functions $q_1$ and $q_2,$ respectively. Then 
		\[D_1 \bar\eta^1J(T_1) \leq \bar \eta^1J(T_2) \leq D_2 \bar \eta^1J(T_1)\]
		where $B_1=\inf_{u\in (0,1)} \left(\frac{\phi(q_2(u))}{\phi(q_1(u)}\right)$ and $B_2=\sup_{u\in (0,1)} \left(\frac{\phi(q_2(u))}{\phi(q_1(u))}\right).$		
	\end{proposition}

	\begin{proposition}
		Let $T_1$ and $T_2$ be the lifetimes of two coherent systems with the same
		structure and with i.d. components having common distributions $F$ and $G$, respectively. If $X \leq_{d} Y$ then we have \[\eta^1J(T_1) \leq \eta^1J(T_2)\]
	\end{proposition}
	
	\begin{corollary}
		Let $T_1$ and $T_2$ be the lifetimes of two coherent systems with the same
		structure and with i.d. components having common distributions $F$ and $G$, respectively. If $X \leq_{hr} Y$ and X or Y is DFR, then \[\eta^1J(T_1) \leq \eta^1J(T_2)\]
	\end{corollary}
	
	\section{Conclusion} \label{s10conclusion}
	Extropy has been studied and generalized by various researchers as discussed in the introduction section. We defined and studied general weighted cumulative extropy, general weighted cumulative residual extropy, general weighted negative cumulative past extropy and conditional weighted past extropy. Non-Parametric estimators for general weighted cumulative extropy, general weighted cumulative residual extropy and weighted negative cumulative extropy have been proposed. We obtained the upper and lower bound of weighted negative cumulative past extropy. We proposed some results to weighted negative cumulative extropy of a coherent system.\\

	\noindent \textbf{\Large Funding} \\
	\\
	Santosh Kumar Chaudhary would like to thank the council of scientific and industrial research (CSIR), Government of India (File Number 09/0081 (14002)/2022-EMR-I) for financial assistance.\\
	\\
	\textbf{ \Large Conflict of interest} \\
	\\
	The authors declare no conflict of interest.


\begin{thebibliography}{99}
		
		\bibitem{lad2015} 
		Lad F., Sanfilippo G., Agro G. Extropy: complementary dual of entropy Statist. Sci., 30 (2015), pp. 40-58
		
		\bibitem{shanon1948} 
		Shannon C.E. A mathematical theory of communication Bell Syst. Tech. J., 27 (1948), pp. 379-423
		
		\bibitem{sudheeshginiindex2021}
		Sudheesh K.K., Dewan I., Sreelaksmi N.
		Non-parametric estimation of gini index with right censored observations
		Statist. Probab. Lett., 175 (2021)
		
		\bibitem{sudheeshetal2022CRE}  
		Sudheesh, K. K. and Sreedevi, E. P. (2022). Non-parametric estimation of cumulative (residual) extropy. Statistics and Probability Letters, 185, 109434.
		
		\bibitem{a}
		\noindent Ash, R., 1990. Information Theory. Dover Publications Inc., New York.
		\bibitem{b}
		\noindent Bagnoli, M. and Bergstrom, T., 2005. Log-concave probability and its applications. Economic Theory, 26, 455-469.
		\bibitem{vasicek} 
		Vasicek, O. ( 1976), A test for normality based on sample entropy. J R Stat Soc Ser B Methodol.38:54–59. 
		\bibitem{tahmasebitoomaj2020}
		Tahmasebi, S.; Toomaj, A. On negative cumulative extropy with applications. Commun. Stat. Theory Methods 2020, 51, 5025–5047.
		\bibitem{balakrishananetal2022}
		Balakrishnan, N.; Buono, F.; Longobardi, M. On Tsallis extropy with an application to pattern recognition. Stat. Probab. Lett. 2022,
		180, 109241. 
		\bibitem{tahmasebietal2022}
		Tahmasebi, S.; Kazemi, M.R.; Keshavarz, A.; Jafari, A.A.; Buono, F. Compressive Sensing Using Extropy Measures of Ranked Set Sampling. Math. Slovaca 2022, accepted for publication.
		\bibitem{kazemietal2021}
		Kazemi, M.R.; Tahmasebi, S.; Buono, F.; Longobardi, M. Fractional Deng Entropy and Extropy and Some Applications. Entropy 2021, 23, 623.
		\bibitem{navarro2022}
		Navarro J (2022) Introduction to system reliability theory. Springer
		\bibitem{barloproschan1975}
		Barlow, R. E., Proschan, F. (1975). Statistical theory of reliability and life testing: probability models. Florida State Univ Tallahassee.
		\bibitem{kazemietal.2022} Kazemi, M.R.; Hashempour, M.; Longobardi, M. Weighted Cumulative Past Extropy and Its Inference. Entropy 2022, 24, 1444. https://doi.org/10.3390/e24101444
		\bibitem{Kattumanniletal.2022} Kattumannil, Sudheesh K. \& E.P., Sreedevi, 2022. "Non-parametric estimation of cumulative (residual) extropy," Statistics \& Probability Letters, Elsevier, vol. 185(C).
		\bibitem{dattaetal.2022} DATTA, S., BANDYOPADHYAY, D., \& SATTEN, G. A. (2010). Inverse Probability of Censoring Weighted U-statistics for Right-Censored Data with an Application to Testing Hypotheses. Scandinavian Journal of Statistics, 37(4), 680–700.
		\bibitem{raoetal2004} Rao, M., Y. Chen, B. C. Vemuri, and F. Wang. 2004. Cumulative residual entropy: A new measure of information. IEEE Transactions on Information Theory 50 (6) :1220–28. doi:10.1109/TIT.2004.
		828057.
		\bibitem{bansalgupta2020} Shilpa Bansal \& Nitin Gupta (2020) Weighted extropies and past extropy of order statistics and k-record values, Communications in Statistics - Theory and Methods, 51:17, 6091-6108, DOI: 10.1080/03610926.2020.1853773
		\bibitem{balaetal2020we} Narayanaswamy Balakrishnan, Francesco Buono \& Maria Longobardi (2022) On weighted extropies, Communications in Statistics - Theory and Methods, 51:18, 6250-6267, DOI: 10.1080/03610926.2020.1860222
		
		\bibitem{guptaskc22} Gupta, N., \& Chaudhary, S. K. (2022). On General Weighted Extropy of Ranked Set Sampling. arXiv preprint arXiv:2207.02003.
		
		\bibitem{janansahietal20}	Jahanshahi, S., Zarei, H., \& Khammar, A. (2020). On cumulative residual extropy. Probability in the Engineering and Informational Sciences, 34(4), 605-625.
		
		\bibitem{shakedshanti2007} Shaked, M., and J. G. Shanthikumar. 2007. Stochastic orders and their applications. San Diego: Academic Press.
		\bibitem{kundunanda2010} Chanchal Kundu \& Asok K. Nanda (2010) Some Reliability Properties of the Inactivity Time, Communications in Statistics - Theory and Methods, 39:5, 899-911, DOI: 10.1080/03610920902807895
		
		\bibitem{navaroetal2013} Navarro, J., Y. del Aguila, M. A. Sordo, and A. Suárez-Llorens. 2013. Stochastic ordering properties for systems with dependent identically distributed components. Applied Stochastic Models in Business and Industry 29 (3):264–78.
		
		\bibitem{kattumani22npecre} Sudheesh, K. K. and Sreedevi, E. P. (2022). Non-parametric estimation of cumulative (residual) extropy. Statistics \& Probability Letters, 185, 109434.
		
		\bibitem{hashempouretal22} M. Hashempour, M. R. Kazemi \& S. Tahmasebi (2022) On weighted cumulative residual extropy: characterization, estimation and testing, Statistics, 56:3, 681-698, DOI: 10.1080/02331888.2022.2072505
		
		\bibitem{qiu2017} Guoxin Qiu, The extropy of order statistics and record values, Statistics \& Probability Letters, Volume 120, 2017, Pages 52-60.
		
		\bibitem{qiujia2018a} Guoxin Qiu, Kai Jia,
		The residual extropy of order statistics,
		Statistics \& Probability Letters, Volume 133,
		2018, Pages 15-22.
		\bibitem{josesathar19} Jitto Jose, E.I. Abdul Sathar,
		Residual extropy of k-record values, Statistics \& Probability Letters, Volume 146, 2019, Pages 1-6.
		
		\bibitem{qiujia2018b} Qiu, G., \& Jia, K. (2018). The residual extropy of order statistics. Statistics \& Probability Letters, 133, 15-22.
		
		\bibitem{qiuetal2019} Qiu, G., Wang, L., \& Wang, X. (2019). ON EXTROPY PROPERTIES OF MIXED SYSTEMS. Probability in the Engineering and Informational Sciences, 33(3), 471-486. doi:10.1017/S0269964818000244
		
		\bibitem{lee2019} Lee, A. J. (2019). U-statistics: Theory and Practice. Routledge.
		
		\bibitem{satharnair21a} Sathar, E.A., Nair, R.D., 2021a. On dynamic weighted extropy. J. Comput. Appl. Math. 393, 113507.
		
		\bibitem{satharnair21b} Sathar, E.A., Nair, R.D., 2021b. On dynamic survival extropy. Comm. Statist. Theory Methods 50, 1295–1313.
		
		\bibitem{satharnair21c} Sathar, E.A., Nair, R.D., 2021c. A study on weighted dynamic survival and failure extropies. Comm. Statist. Theory Methods 1–20.
		
		\bibitem{kattumannilsudheeshetal21} Sudheesh, K.K., Dewan, I., Sreelaksmi, N., 2021. Non-parametric estimation of gini index with right censored observations. Statist. Probab. Lett. 175, 109113. http://dx.doi.org/10.1016/j.spl.2021.109113.
		
		\bibitem{dataaetal2010} Datta S, Bandyopadhyay D, Satten GA. Inverse probability of censoring weighted U-statistics for right-censored data with an application to testing hypotheses. Scandinavian J Statist.	2010;37:680–700.
		
		\bibitem{lehman1951} E. L. Lehmann. "Consistency and Unbiasedness of Certain Nonparametric Tests." Ann. Math. Statist. 22 (2) 165 - 179, June, 1951. https://doi.org/10.1214/aoms/1177729639
	\end{thebibliography}
\end{document}